\documentclass[12pt]{article}

\usepackage{amssymb}
\usepackage{latexsym}

\newtheorem{definition}{Definition}[section]
\newtheorem{theorem}[definition]{Theorem}
\newtheorem{lemma}[definition]{Lemma}
\newtheorem{corollary}[definition]{Corollary}
\newtheorem{remark}[definition]{Remark}

\newtheorem{conjecture}[definition]{Conjecture}
\newtheorem{problem}[definition]{Problem}
\newtheorem{note}[definition]{Note}

\newtheorem{proposition}[definition]{Proposition}

\def\I{\mathbb I}

\def\F{\mathbb F}
\def\K{\mathbb F}

\begin{document}

\title{\bf Tridiagonal pairs of   \\
Krawtchouk type
}
\author{
Tatsuro Ito{\footnote{
Department of Computational Science,
Faculty of Science,
Kanazawa University,
Kakuma-machi,
Kanazawa 920-1192, Japan
}}
$\;$ and
Paul Terwilliger{\footnote{
Department of Mathematics, University of
Wisconsin, 480 Lincoln Drive, Madison WI 53706-1388 USA}
}}
\date{}

\maketitle
\begin{abstract}
Let $\K$ denote an algebraically closed field with
characteristic 0 and let $V$ denote a vector space over
$\K$ with finite positive dimension. Let $A,A^*$ 
denote a tridiagonal pair on $V$ with diameter $d$.
We say that $A,A^*$ has {\it Krawtchouk type} whenever
the sequence $\lbrace d-2i\rbrace_{i=0}^d$
is a standard
ordering of the eigenvalues of $A$ and a standard ordering
of the eigenvalues of $A^*$. Assume $A,A^*$ has
Krawtchouk type.
We show that there exists a nondegenerate
symmetric bilinear form $\langle \,,\,\rangle $ on $V$ such that
$\langle Au,v\rangle= \langle u,Av \rangle$ and
$\langle A^*u,v \rangle= \langle u,A^*v \rangle$ 
for $u,v\in V$.
We show that the following tridiagonal pairs are isomorphic:
(i) $A,A^*$;
(ii) $-A,-A^*$;
(iii) $A^*,A$;
(iv) $-A^*,-A$.
We give a number of related results and conjectures.

\bigskip
\noindent
{\bf Keywords}. 
Tridiagonal pair, Leonard pair, tetrahedron Lie algebra.
 \hfil\break
\noindent {\bf 2000 Mathematics Subject Classification}. 
Primary: 33C45. Secondary: 
05E30, 05E35,
15A21, 17B65.
 \end{abstract}

\section{Tridiagonal pairs}

\noindent 
Throughout this paper $\K$ denotes a field.

\medskip
\noindent 
We begin by recalling the notion of a tridiagonal pair. 
We will use the following terms.
Let $V$  denote
a vector space over $\K$ with finite positive dimension.
Let ${\rm End}(V)$ denote the $\K$-algebra of all linear
transformations from $V$ to $V$.
For a subspace $W \subseteq V$ and
$A \in 
{\rm End}(V)$,
we call $W$ an
 {\it eigenspace} of $A$ whenever 
 $W\not=0$ and there exists $\theta \in \K$ such that 
\begin{eqnarray*}
W=\lbrace v \in V \;\vert \;Av = \theta v\rbrace.
\end{eqnarray*}
We say that $A$ is {\it diagonalizable} whenever
$V$ is spanned by the eigenspaces of $A$.

\begin{definition}  
{\rm \cite{TD00}}
\label{def:tdp}
\rm
Let $V$ denote a vector space over $\K$
with finite positive dimension.
By a {\it tridiagonal pair} (or {\it TD pair}) on $V$
we mean an ordered pair $A,A^*$ of elements in
${\rm End}(V)$
that satisfy the following four conditions.
\begin{enumerate}
\item Each of $A,A^*$ is diagonalizable.
\item There exists an ordering $\lbrace V_i\rbrace_{i=0}^d$ of the  
eigenspaces of $A$ such that 
\begin{equation}
A^* V_i \subseteq V_{i-1} + V_i+ V_{i+1} \qquad \qquad (0 \leq i \leq d),
\label{eq:t1}
\end{equation}
where $V_{-1} = 0$ and $V_{d+1}= 0$.
\item There exists an ordering $\lbrace V^*_i\rbrace_{i=0}^{\delta}$ of
the  
eigenspaces of $A^*$ such that 
\begin{equation}
A V^*_i \subseteq V^*_{i-1} + V^*_i+ V^*_{i+1} 
\qquad \qquad (0 \leq i \leq \delta),
\label{eq:t2}
\end{equation}
where $V^*_{-1} = 0$ and $V^*_{\delta+1}= 0$.
\item There does not exist a subspace $W$ of $V$ such  that $AW\subseteq W$,
$A^*W\subseteq W$, $W\not=0$, $W\not=V$.
\end{enumerate}
We say the pair $A,A^*$ is {\it over $\K$}.
 We call $V$ the {\it vector space 
underlying $A,A^*$}.
We call ${\rm End}(V)$ the {\it ambient algebra}
of $A,A^*$.
\end{definition}

\begin{note}
\rm
According to a common notational convention $A^*$ denotes 
the conjugate-transpose of $A$. We are not using this convention.
In a TD pair $A,A^*$ the linear transformations $A$ and $A^*$
are arbitrary subject to (i)--(iv) above.
\end{note}

\noindent We refer the reader to 
\cite{hasan},
\cite{hasan2},
\cite{TD00},
\cite{shape},
\cite{tdanduq},
\cite{NN},
\cite{N:refine},
\cite{nomsplit},
\cite{madrid} and the references therein 
for background on tridiagonal pairs. 

\medskip
\noindent In order to motivate
our results we recall a few facts about TD pairs.
Let $A,A^*$ denote a TD pair on $V$, as in Definition
\ref{def:tdp}.
It turns out that the integers $d$ and $\delta$ from
conditions (ii) and (iii) are equal
\cite[Lemma~4.5]{TD00}; we call this common value
the {\it diameter} of $A,A^*$.
An ordering of the eigenspaces of $A$ (resp. $A^*$)
will be called
{\it standard} whenever it satisfies
(\ref{eq:t1}) (resp. 
(\ref{eq:t2})). 
We comment on the uniqueness of the standard ordering.
Let 
 $\lbrace V_i\rbrace_{i=0}^d$ denote a standard ordering
of the eigenspaces of $A$.
Then the ordering 
 $\lbrace V_{d-i}\rbrace_{i=0}^d$ 
is standard and no other ordering is standard.
A similar result holds for the 
 eigenspaces of $A^*$.
An ordering of the eigenvalues of $A$ (resp. $A^*$)
will be called {\it standard} whenever the 
corresponding ordering of the eigenspaces of
$A$ (resp. $A^*$) is standard.
Let
 $\lbrace V_i\rbrace_{i=0}^d$
 (resp. $\lbrace V^*_i\rbrace_{i=0}^d$)
denote a standard ordering
of the eigenspaces of $A$ (resp. $A^*$).
By \cite[Corollary~5.7]{TD00},  for $0 \leq i \leq d$ 
the spaces $V_i$,  $V^*_i$ have the same dimension;
we denote this common dimension by $\rho_i$.
By the construction $\rho_i\not=0$.
 By  \cite[Corollary~5.7]{TD00} 
and \cite[Corollary~6.6]{TD00}
the sequence $\lbrace \rho_i\rbrace_{i=0}^d$
is symmetric and unimodal; that is
$ \rho_i =\rho_{d-i}$ for $0 \leq i \leq d$
and 
$\rho_{i-1} \leq \rho_{i}$ for  $1 \leq i \leq d/2$.
We call the sequence
$\lbrace \rho_i\rbrace_{i=0}^d$
  the {\it shape} of $A,A^*$.

\medskip
\noindent The following special case has received a lot of
attention. By a {\it Leonard pair} we mean
a TD pair with shape
$(1,1,\ldots, 1)$
\cite[Definition~1.1]{LS99},
  \cite[Lemma~2.2]{qSerre}.
There is a natural correspondence between
the Leonard pairs and 
a family of orthogonal
polynomials consisting of the $q$-Racah polynomials
\cite{AWil},
\cite{GR} and their relatives
\cite{KoeSwa}.
This family coincides
with the terminating branch of the Askey scheme
\cite{TLT:array}. 
See
\cite{N:aw},
\cite{NT:balanced},
\cite{NT:formula},
\cite{NT:det},
\cite{NT:mu},
\cite{NT:span},
\cite{NT:switch},
\cite{nomsplit},
\cite{qSerre},
\cite{LS24},
\cite{conform},
\cite{lsint},
\cite{TLT:split},
	  \cite{qrac},
	  \cite{aw}
	   for more information about Leonard pairs.

\medskip
\noindent We mention some notation for later use.
For an indeterminate $\lambda$ let $\K \lbrack \lambda \rbrack$
denote the $\K$-algebra of all polynomials in $\lambda$ that
have coefficients in $\K$.

\medskip
\noindent
We now summarize our results. For the rest
of this section assume $\K$ is algebraically closed with
characteristic 0. Let $A,A^*$ denote a TD pair over
$\K$ and let $d$ denote the diameter. We say that
$A, A^*$ has {\it Krawtchouk type}
whenever 
the sequence 
$\lbrace d-2i\rbrace_{i=0}^d$
is a standard ordering of the eigenvalues of $A$
and a
standard ordering of the eigenvalues of $A^*$.
Assume $A,A^*$ has Krawtchouk type and let $V$ denote the
underlying vector space. We will prove the following.
\begin{itemize}
\item
Let $\lbrace \rho_i\rbrace_{i=0}^d$ denote the shape
of $A,A^*$. 
 Then there exists a nonnegative integer $N$ and
  positive integers $d_1, d_2, \ldots, d_N$ such that
\begin{eqnarray*}
  \sum_{i=0}^d \rho_i \lambda^i
  = \prod_{j=1}^N
   (1+\lambda+\lambda^2+\cdots + \lambda^{d_j}).
   \end{eqnarray*}
\item
There exists a nonzero bilinear form
 $\langle\,,\,\rangle$ on $V$ such that
$\langle Au,v\rangle= \langle u,Av \rangle $
and 
$\langle A^*u,v \rangle= \langle u,A^*v \rangle $
for $u,v\in V$. The form is unique up to multiplication
by a nonzero scalar in $\K$. The form is nondegenerate and
symmetric.
\item
There exists a unique antiautomorphism $\dagger $ 
of ${\rm End}(V)$
that fixes each of $A,A^*$. Moreover $X^{\dagger \dagger}=X $
for all $X \in {\rm End}(V)$.
\item
Let $B$ (resp. $B^*$) denote the image
of $A$ (resp. $A^*$) under the canonical 
 anti-isomorphism
${\rm End}(V)\to 
{\rm End}({\tilde V})$, where $\tilde V$ is
the vector space dual to $V$.
Then the TD pairs
$A,A^*$ and 
 $B,B^*$ are isomorphic.
\item The following TD pairs are mutually isomorphic:
\begin{eqnarray*}
A,A^*; \qquad \qquad 
-A,-A^*; \qquad \qquad 
A^*,A; \qquad \qquad 
-A^*,-A.
\end{eqnarray*}
\end{itemize}
In addition we 
prove the following.
We associate with $A,A^*$ a polynomial
$P_{A,A^*}$ in $\K\lbrack \lambda \rbrack$
called the Drinfel'd polynomial. We show that
the map $A,A^*\mapsto P_{A,A^*}$ induces a bijection
between 
(i) the set of isomorphism classes of TD pairs over $\K$
that have Krawtchouk type;
(ii) the set of polynomials in $\K\lbrack \lambda \rbrack$ 
that have constant coefficient
1 and do not vanish at $\lambda = 1$.

\medskip
\noindent 
To obtain the above results 
we use a bijection
due to Hartwig
\cite{Ha}
involving the TD pairs of Krawtchouk type and the finite-dimensional
irreducible modules for the tetrahedron algebra
\cite{HT}.

\medskip
\noindent We finish the paper with a number of conjectures and
problems concerning general TD pairs.

\medskip
\noindent The paper is organized as follows.
In Sections 2--5 we review some facts about general
TD pairs.
In Sections 6, 7 we recall the tetrahedron algebra and
its relationship to 
the TD pairs of Krawtchouk type. We use this relationship
in Sections 
8--13 to prove our results on TD pairs of Krawtchouk type.
In Section 14 we give our conjectures and open problems.

\medskip
\noindent 
We remark that Alnajjar and Curtin
\cite{CurtH} obtained some results similar to
ours for the
 TD pairs of $q$-Serre type.

\section{Isomorphisms of TD pairs}

\noindent In this section we consider 
the concept of isomorphism
for TD pairs from several points of view.

\begin{definition}
\label{def:iso}
\rm
Let $A,A^*$ and $B,B^*$ denote TD pairs over $\K$.
By an {\it isomorphism of TD pairs}  from
$A,A^*$ to $B,B^*$ we mean
a vector space isomorphism $\gamma$ from
the vector space underlying $A,A^*$ to
the vector space underlying $B,B^*$ such that
both
\begin{eqnarray*}
\label{eq:iso}
\gamma A
=B\gamma,
\qquad \qquad 
\gamma A^*=
B^*\gamma.
\end{eqnarray*}
\end{definition}

\noindent We have a comment.

\begin{lemma}
\label{lem:id}
Assume $\K$ is algebraically closed.
Let $V$ denote a vector space over $\K$ with finite
positive dimension and let $A,A^*$ denote a TD pair
on $V$. Then the following (i), (ii) are equivalent
for 
$\gamma \in 
{\rm End}(V)$.
\begin{enumerate}
\item $\gamma$ commutes with each of $A, A^*$;
\item there 
exists $\alpha \in \K$ such that
$\gamma =\alpha I $.
\end{enumerate}
\end{lemma}
\noindent {\it Proof:}
$(i) \Rightarrow (ii)$
Since $V$ has finite positive dimension
and since 
 $\K$ is algebraically closed
there exists at least one eigenspace $W\subseteq V$
for 
$\gamma$.  
Note that
$AW\subseteq W$ since
$\gamma A=A \gamma$ and
$A^*W\subseteq W$ since
$\gamma A^*=A^* \gamma$.
 Also $W\not=0$ by the definition of
an eigenspace so $W=V$  
 in view of
Definition \ref{def:tdp}(iv). The result follows.

\noindent 
$(ii) \Rightarrow (i)$ Clear.
\hfill $\Box $ \\

\begin{proposition}
\label{prop:id}
Assume $\K$ is algebraically closed.
Let $V$ denote a vector space over $\K$ with finite
positive dimension and let $A,A^*$ denote a TD pair
on $V$. Then the following (i), (ii) are equivalent
for 
$\gamma \in 
{\rm End}(V)$.
\begin{enumerate}
\item
$\gamma$ is an 
isomorphism
of TD pairs from $A,A^*$ to $A,A^*$;
\item there 
exists a nonzero  $\alpha \in \K$ such that
$\gamma =\alpha I $.
\end{enumerate}
\end{proposition}
\noindent {\it Proof:}
$(i) \Rightarrow (ii)$
The map $\gamma$ is nonzero and commutes with
each of $A,A^*$ so the result follows in view of
Lemma \ref{lem:id}.

\noindent 
$(ii) \Rightarrow (i)$ Clear.
\hfill $\Box $ \\ 

\noindent 
We have been discussing the concept of isomorphism for TD pairs.
In this discussion we now shift the emphasis 
from maps involving the underlying vector spaces
to maps involving the ambient algebras.
We do this as follows.
Let $V$ denote a vector space over $\K$ with finite
positive dimension and let $\gamma : V\to V'$ denote
an isomorphism of vector spaces. Observe that
there exists an $\K$-algebra isomorphism $\varrho:
{\rm End}(V)\to 
{\rm End}(V')$
such that $X^{\varrho}= \gamma X \gamma^{-1}$ for
all $X \in
{\rm End}(V)$.
Conversely let 
$\varrho:  
{\rm End}(V)\to 
{\rm End}(V')$ denote an $\K$-algebra isomorphism.
Then by the
Skolem-Noether theorem \cite[Corollary~9.122]{rotman}
there exists an isomorphism of vector spaces
$\gamma:V \to V'$ such that
$X^{\varrho}= \gamma X \gamma^{-1}$ for
all $X \in
{\rm End}(V)$. Combining these comments with
Definition \ref{def:iso} we obtain the following result.

\begin{corollary}
\label{cor:iso2}
Let $A,A^*$ and $B,B^*$ denote TD pairs over $\F$.
Then these TD pairs are isomorphic if and only if
there exists an $\F$-algebra isomorphism from
the ambient algebra of $A,A^*$ to the
the ambient algebra of $B,B^*$ that sends
$A$ to $B$ and
$A^*$ to $B^*$.
\end{corollary}

\noindent Let $V$ denote a vector space over $\F$ with
finite positive dimension. By an {\it automorphism}
of ${\rm End}(V)$ we mean an $\K$-algebra isomorphism
from 
${\rm End}(V)$ to 
${\rm End}(V)$.

\begin{corollary}
\label{cor:id2}
Assume $\F$ is algebraically closed.
Let $V$ denote a vector space over $\K$ with finite
positive dimension and let $A,A^*$ denote a TD pair
on $V$. Then the identity map is the unique
automorphism of
${\rm End}(V)$ that fixes each of $A, A^*$.
\end{corollary}
\noindent {\it Proof:} 
Immediate from
Proposition
\ref{prop:id} and our comments above
Corollary
\ref{cor:iso2}.
\hfill $\Box $ \\ 

\section{Bilinear forms}

\noindent Later in the paper we will discuss some bilinear
forms related to TD pairs; to prepare for this we recall some
relevant linear algebra.

\medskip
\noindent
Let $V$ and $V'$ denote 
vector spaces over $\K$ with the same finite positive dimension.
A  map $\langle \,,\,\rangle : V \times V' \rightarrow \K$
is called a 
{\it bilinear form}
whenever 
the following four conditions hold for
$u,v \in V$, for
$u',v' \in V'$,   
and for $\alpha \in \K$:
(i) $\langle u+v,u' \rangle =
 \langle u,u' \rangle
 +
 \langle v,u'\rangle$;
 (ii)
 $
 \langle \alpha u,u' \rangle
 =
 \alpha \langle u,u' \rangle
 $;
 (iii)
  $\langle u,u'+v' \rangle =
   \langle u,u' \rangle
   +
   \langle u,v' \rangle$;
   (iv)
   $
   \langle u, \alpha u' \rangle
   =
   \alpha \langle u,u' \rangle
   $.
   We observe
   that a scalar multiple of
   a bilinear form is a bilinear form.
Let
   $
   \langle\, ,\, \rangle :
   V \times V' \to \K
   $
   denote a bilinear form.
   Then the following
   are equivalent: (i) there exists a nonzero $v \in V$ such
   that $
   \langle v,v' \rangle
   = 0
   $
   for all $v' \in V'$;
    (ii) there exists a nonzero $v' \in V'$ such
    that $
    \langle v,v' \rangle
    = 0
    $
    for all $v \in V$.
    The form
    $\langle \,,\,\rangle $
    is
    said to be {\it degenerate }
    whenever (i), (ii) hold and {\it nondegenerate}
    otherwise.
Here is an example of a nondegenerate bilinear form.
Let ${\tilde V}$ denote the vector space dual
to $V$, consisting
of the linear transformations from
$V$ to $\K$. 
Define a map
 $\langle \,,\,\rangle : V \times {\tilde V} \rightarrow \K$
such that
 $\langle v,f\rangle =f(v)$
for all $v \in V$ and $f \in {\tilde V}$.
Then 
 $\langle \,,\,\rangle$ is a nondegenerate bilinear form.
We call this form the {\it canonical} bilinear form
between $V$ and ${\tilde V}$.
By a {\it bilinear form on $V$} we mean
a bilinear form
   $
   \langle\, ,\, \rangle :V\times V\to \F$.
   This form is said to be {\it symmetric}
whenever
   $
   \langle u,v \rangle
   =
   \langle v,u \rangle
   $
   for all $u, v \in V$.

\medskip
\noindent 
Again let $V$ and $V'$ denote  vector spaces over $\K$ with
the same finite positive dimension.
By an {\it $\K$-algebra
\it anti-isomorphism} from
${\rm{End}}(V)$ to 
${\rm{End}}(V')$ we mean 
an $\K$-linear bijection $\varrho :
{\rm{End}}(V) \to
{\rm{End}}(V')$
such that
$(XY)^{\varrho}= Y^\varrho X^\varrho$
for all $X,Y \in 
{\rm{End}}(V)$.
Let 
 $\langle \,,\,\rangle : V \times V' \rightarrow \K$
denote a nondegenerate bilinear form.
Then there exists
a unique $\K$-algebra
anti-isomorphism $\varrho:
{\rm{End}}(V) \to 
{\rm{End}}(V')$ 
such that
$\langle Xv,v'\rangle = 
\langle v,X^\varrho u'\rangle $
for all $X \in {\rm{End}}(V)$, $v \in V$, $v' \in V'$.
We say that $\varrho$  {\it corresponds}
to $\langle \,,\,\rangle $.
By the {\it canonical} $\K$-algebra
anti-isomorphism
from ${\rm{End}}(V)$ to 
${\rm{End}}({\tilde V})$
we mean the one that
corresponds
to the canonical bilinear form.
By an 
{\it antiautomorphism} of 
${\rm{End}}(V)$ we mean an 
$\K$-algebra anti-isomorphism from 
${\rm{End}}(V)$ to
${\rm{End}}(V)$.

\section{New TD pairs from old}

\noindent 
Let $V$ denote a vector space over $\K$ with
finite positive dimension and let
 $A,A^*$ denote a TD pair on $V$.
We can modify this pair in several ways to get another
TD pair. For instance
the ordered pair $A^*,A$ is a TD pair on $V$ which
is potentially nonisomorphic to $A,A^*$.
Let $\alpha, \alpha^*, \beta, \beta^*$ denote scalars in
$\K$
such that each of $\alpha, \alpha^*$ is nonzero.
Then the pair
\begin{eqnarray*}
\alpha A + \beta I, \qquad \qquad  
\alpha^* A^* + \beta^* I
\end{eqnarray*}
is a TD pair on $V$ which is potentially
nonisomorphic to $A,A^*$. Let $\tilde V$ denote the
vector space dual to $V$. 
 Let
$B$ (resp. $B^*$) denote the image of
$A$ (resp. $A^*$) under the canonical 
$\K$-algebra anti-isomorphism
${\rm End}(V)\to 
{\rm End}({\tilde V})$
from Section 3.
By \cite[Theorem~1.2]{CurtH} 
the pair $B,B^*$ is
a TD pair on $\tilde V$ which is potentially nonisomorphic
to $A,A^*$.

\section{The split decomposition of a TD pair} 

\medskip
\noindent 
Let $V$ denote a vector space over $\K$ with finite
positive dimension and
let $A,A^*$ denote a TD pair on $V$. 
Let $\lbrace \theta_i\rbrace_{i=0}^d$
(resp. 
 $\lbrace \theta^*_i\rbrace_{i=0}^d$)
denote a standard ordering of the eigenvalues
for $A$ (resp. $A^*$). We recall the corresponding
split decomposition and split sequence.
Let $\lbrace V_i\rbrace_{i=0}^d$
(resp. $\lbrace V^*_i\rbrace_{i=0}^d$)
denote the 
ordering
of the eigenspaces of $A$ (resp. $A^*$) associated with
$\lbrace \theta_i\rbrace_{i=0}^d$
(resp. 
 $\lbrace \theta^*_i\rbrace_{i=0}^d$).
For $0 \leq i \leq d$ define
\begin{eqnarray*}
\label{eq:uidef}
U_i = (V^*_0+V^*_1+\cdots + V^*_i)\cap
(V_i+V_{i+1}+\cdots + V_d).
\end{eqnarray*}
By \cite[Theorem~4.6]{TD00}
\begin{eqnarray*}
V = \sum_{i=0}^d U_i   \qquad \qquad (\mbox{\rm direct sum}),
\end{eqnarray*}
and for $0 \leq i \leq d$ both
\begin{eqnarray*}
U_0+U_1+\cdots +U_i &=&V^*_0+V^*_1+\cdots + V^*_i,\\
U_i+U_{i+1}+\cdots +U_d &=&V_i+V_{i+1}+\cdots + V_d.
\end{eqnarray*}
By \cite[Corollary~5.7]{TD00}
$U_i$  has dimension $\rho_i$ where
$\lbrace \rho_i\rbrace_{i=0}^d$ is the shape of $A,A^*$.
By \cite[Theorem~4.6]{TD00} both
\begin{eqnarray}
\label{eq:ath}
(A-\theta_iI)U_i &\subseteq & U_{i+1}, \\
\label{eq:asths}
(A^*-\theta^*_iI)U_i & \subseteq & U_{i-1},
\end{eqnarray}
where $U_{-1}=0$ and $U_{d+1}=0$.
The sequence $\lbrace U_i\rbrace_{i=0}^d$
is called the {\it split decomposition} of $V$
with respect to 
$\lbrace \theta_i\rbrace_{i=0}^d$
and $\lbrace \theta^*_i\rbrace_{i=0}^d$ \cite[Section~4]{TD00}.
Now assume 
$\rho_0=1$, so that $U_0$ has dimension 1.
For $0 \leq i \leq d$ the space $U_0$ is invariant
under
\begin{eqnarray}
\label{eq:aaseig}
(A^*-\theta^*_1I) 
(A^*-\theta^*_2I)  \cdots 
(A^*-\theta^*_iI) 
(A-\theta_{i-1}I) \cdots 
(A-\theta_1I)   
(A-\theta_0I);
\end{eqnarray}
let $\zeta_i$ denote the corresponding eigenvalue.
We call the sequence $\lbrace \zeta_i\rbrace_{i=0}^d$
the {\it split sequence} of $A,A^*$
with respect to 
$\lbrace \theta_i\rbrace_{i=0}^d$
and $\lbrace \theta^*_i\rbrace_{i=0}^d$.

\begin{note}
\rm
In the literature on Leonard pairs there
are two sequences of scalars
called the first split sequence and the
second split sequence \cite[Section~3]{LS99}. 
These sequences are related to
the above split sequence as follows. 
Let $A,A^*$ denote a Leonard pair and
fix a standard ordering $\lbrace \theta_i \rbrace_{i=0}^d$
(resp.
 $\lbrace \theta^*_i \rbrace_{i=0}^d$) of the eigenvalues
 of $A$ (resp. $A^*$).
Let $\lbrace \varphi_i\rbrace_{i=1}^d$ 
(resp. 
 $\lbrace \phi_i\rbrace_{i=1}^d$) 
denote the
corresponding first split sequence (resp. second split sequence)
in the sense of \cite{LS99}.
Then the sequence
 $\lbrace \varphi_1\varphi_2\cdots \varphi_i\rbrace_{i=0}^d$
(resp. 
 $\lbrace \phi_1\phi_2\cdots \phi_i\rbrace_{i=0}^d$)
is the split sequence of $A,A^*$ associated with
$\lbrace \theta_i \rbrace_{i=0}^d$ and
 $\lbrace \theta^*_i \rbrace_{i=0}^d$
(resp. 
$\lbrace \theta_{d-i} \rbrace_{i=0}^d$ and
 $\lbrace \theta^*_i \rbrace_{i=0}^d$).
\end{note}


\section{The tetrahedron algebra}

\noindent From now until the end of Section 13 we assume:

\medskip
\centerline{The field $\K$ is algebraically closed with characteristic 0.}

\medskip
\noindent 
We now recall some facts about the tetrahedron
 algebra that we will use later in the paper.
We start with a definition.

\begin{definition}
\label{def:tet}
\rm
\cite[Definition~1.1]{HT}
Let $\boxtimes$ denote the Lie algebra over $\K$
that has  generators
\begin{eqnarray}
\label{eq:boxgen}
\lbrace x_{ij} \,|\,i,j\in \I, i\not=j\rbrace
\qquad \qquad \I = \lbrace 0,1,2,3\rbrace
\end{eqnarray}
and the following relations:
\begin{enumerate}
\item[{\rm (i)}]  
For distinct $i,j\in \I$,
\begin{eqnarray*}
x_{ij}+x_{ji} = 0.
\label{eq:rel0}
\end{eqnarray*}
\item[{\rm (ii)}]  
For mutually distinct $i,j,k\in \I$,
\begin{eqnarray*}
\lbrack x_{ij},x_{jk}\rbrack = 2x_{ij}+2x_{jk}.
\label{eq:rel1}
\end{eqnarray*}
\item[{\rm (iii)}]  
For mutually distinct $i,j,k,\ell \in \I$,
\begin{eqnarray}
\label{eq:dg}
\lbrack x_{ij},
\lbrack x_{ij},
\lbrack x_{ij},
x_{k \ell}\rbrack \rbrack \rbrack= 
4 \lbrack x_{ij},
x_{k \ell}\rbrack.
\label{eq:rel2}
\end{eqnarray}
\end{enumerate}
We call $\boxtimes$ the {\it tetrahedron algebra}
or ``tet'' for short. 
\end{definition}

\noindent 
Let $V$ denote a finite-dimensional
irreducible $\boxtimes$-module.
We recall the diameter of $V$  and the shape of $V$.
By 
\cite[Theorem~3.8]{Ha}
each generator $x_{ij}$ of $\boxtimes$ is diagonalizable   
on $V$.
Also by \cite[Theorem~3.8]{Ha}
there exists an integer $d\geq 0$ such that
for each generator $x_{ij}$ the set of distinct eigenvalues
on $V$ is $\lbrace d-2n\,|\,0 \leq n \leq d\rbrace$. 
We call
$d$ the {\it diameter} of $V$.
By \cite[Corollary~3.6]{Ha}, 
for $0 \leq n\leq d$ there exists a positive integer $\rho_n$ such
that for each generator
 $x_{ij}$ 
the $(2n-d)$-eigenspace in $V$ has dimension
$\rho_n$.
We call the sequence 
$\lbrace \rho_n \rbrace_{n=0}^d$  the
{\it shape} of $V$.

\begin{theorem} {\rm \cite[Theorem~16.5]{Ev}}
\label{lem:shape}
 Let 
 $\lbrace \rho_n\rbrace_{n=0}^d$ denote the shape of
a finite-dimensional irreducible
$\boxtimes$-module.
Then there exists a nonnegative integer $N$ and
positive integers $\lbrace d_j\rbrace_{j=1}^N$ such that
\begin{eqnarray*}
\sum_{n=0}^d \rho_n \lambda^n = \prod_{j=1}^N 
(1+\lambda+\lambda^2+\cdots + \lambda^{d_j}).
\end{eqnarray*}
In particular
$\rho_0=1$.
\end{theorem}

\medskip
\noindent We recall how the symmetric group $S_4$
acts on $\boxtimes$ as a group of automorphisms.
We identify $S_4$ with the 
group of permutations of $\I$.
We use the cycle notation; for example
$(1,2,3)$ denotes the element of $S_4$ that
sends
 $1\mapsto 2
 \mapsto 3
 \mapsto 1$ and
$0\mapsto 0$.
Note that $S_4$ acts on the set of generators for
$\boxtimes$ by 
permuting the indices:
\begin{eqnarray*}
\sigma(x_{ij}) = x_{\sigma(i),\sigma(j)}
\qquad \qquad \sigma \in S_4, \qquad i,j\in \I, \quad i\not=j.
\end{eqnarray*}
This action leaves invariant the defining relations for $\boxtimes$
and therefore induces an
 action of $S_4$ on $\boxtimes$ as a group
of automorphisms.

\medskip
\noindent 
Let $G$ denote the unique normal subgroup of
$S_4$ that has cardinality 4. $G$ consists of 
\begin{eqnarray*}
(01)(23),
\qquad \qquad 
(02)(13),
\qquad \qquad 
(03)(12)
\end{eqnarray*}
and the identity element.

\begin{theorem}
\label{thm:bill} 
{\rm \cite[Theorem~19.1]{Ev}}
Let $V$ denote a finite-dimensional irreducible
$\boxtimes$-module. For a nonidentity
$\sigma \in G$ there exists a nonzero bilinear form
$\langle \,,\,\rangle_{\sigma}$ on $V$ such that
\begin{eqnarray}
\langle \xi.u,v\rangle_{\sigma}
=
-\langle u,\sigma(\xi).v\rangle_{\sigma}
\qquad \qquad 
\xi \in \boxtimes, \qquad u,v\in V.
\end{eqnarray}
This form is unique up to multiplication by a nonzero
scalar in $\F$. This form is nondegenerate and
symmetric.
\end{theorem}

\noindent 
Let $V$ denote a $\boxtimes$-module. For $\sigma \in S_4$ 
there exists a $\boxtimes$-module structure on $V$, called
{\it $V$ twisted via $\sigma$}, that behaves as follows:
for all $\xi \in \boxtimes$ and $v \in V$
the vector $\xi.v$ computed in
$V$ twisted via $\sigma$ coincides with the vector
$\sigma^{-1}(\xi).v$ computed in
the original $\boxtimes$-module  $V$.
See \cite[Section~7]{Ev} for more information
on twisting.

\begin{theorem}
\label{thm:G}
{\rm \cite[Corollary~15.5]{Ev}} 
Let $V$ denote a finite-dimensional irreducible
$\boxtimes$-module.
 Then for $\sigma \in G$ the following
are isomorphic:
\begin{enumerate}
\item the $\boxtimes$-module $V$ twisted via $\sigma$;
\item the $\boxtimes$-module $V$.
\end{enumerate}
\end{theorem}

\noindent
Let $V$ denote a finite-dimensional irreducible
$\boxtimes$-module.
We recall the corresponding Drinfel'd polynomial.
Abbreviate 
\begin{eqnarray}
e^+:=\frac{x_{01}+x_{20}}{2},
\qquad \qquad
e^-:=\frac{x_{23}+x_{02}}{2}.
\label{eq:epm}
\end{eqnarray}
Let $U$ denote the eigenspace of $x_{20}$
on $V$ for the eigenvalue $-d$, where $d$ denotes the
diameter of $V$.
Note that $U$ has
dimension 1 by the last line of Theorem
\ref{lem:shape}.
For $0 \leq i \leq d$ the
space $U$ is invariant under
 $(e^-)^i(e^+)^i$
\cite[Lemma~17.1]{Ev}; let
$\vartheta_i=\vartheta_i(V)$ denote the corresponding eigenvalue.
Define a polynomial $P_V \in \K\lbrack \lambda \rbrack$ by 
\begin{eqnarray}
\label{eq:drin}
P_V = \sum_{i=0}^d
\frac{(-1)^i \vartheta_i \lambda^i}{(i!)^2}.
\end{eqnarray}
We call $P_V$ the
{\it Drinfel'd polynomial} of $V$
\cite[Definition~17.3]{Ev}.

\begin{theorem}
\label{thm:drintet}
{\rm \cite[Corollary~17.7]{Ev}}
The map $V \mapsto P_V$ induces a bijection between the following
two sets:
\begin{enumerate}
\item the isomorphism classes of finite-dimensional irreducible
$\boxtimes$-modules;
\item the polynomials in $\K\lbrack \lambda \rbrack$ that
have constant coefficient 1 and are nonzero at $\lambda=1$.
\end{enumerate}
\end{theorem}
See 
 \cite{BT},
 \cite{E},
\cite{Ha}, 
\cite{HT},
\cite{qtet},
\cite{ITdrg},
\cite{Ev},
\cite{PT}
for more background information on $\boxtimes$ and related topics.

\section{TD pairs of Krawtchouk type and the tetrahedron algebra}

\noindent In \cite{Ha} Hartwig relates
the TD pairs of Krawtchouk type to 
the finite-dimensional 
irreducible $\boxtimes$-modules.
His results are summarized in the following two theorems 
and subsequent remark.

\begin{theorem}
\label{prop:ha1}
{\rm \cite[Theorem~1.7]{Ha}}
Let $V$ denote a finite-dimensional irreducible $\boxtimes$-module.
Then the generators $x_{01}, x_{23}$ act on $V$ as a 
TD pair of Krawtchouk type.
\end{theorem}

\begin{theorem}
\label{prop:ha2}
{\rm \cite[Theorem~1.8]{Ha}}
Let $V$ denote a vector space over $\K$ with finite positive
dimension and let $A,A^*$ denote a TD pair on $V$ that has 
Krawtchouk type. Then there exists a unique $\boxtimes$-module
structure on $V$ such that the generators $x_{01}, x_{23}$
act on $V$ as $A,A^*$ respectively.
This $\boxtimes$-module is irreducible.
\end{theorem}

\begin{remark}
\label{prop:ha3}
\rm
\cite[Remark~1.9]{Ha}
Combining the previous two theorems we obtain a bijection
between the following two sets:
\begin{enumerate}
\item the isomorphism classes of TD pairs over $\K$ that have
Krawtchouk type;
\item
the isomorphism classes of finite-dimensional irreducible
$\boxtimes$-modules.
\end{enumerate}
\end{remark}

\section{TD pairs of Krawtchouk type; the shape}

\noindent In this section we describe the shape
of a TD pair that has Krawtchouk type.
We start with two observations.

\begin{lemma}
\label{lem:diam2}
Let $V$ denote a finite-dimensional irreducible
$\boxtimes$-module. Then the following coincide:
\begin{enumerate}
\item the diameter of the $\boxtimes$-module $V$ from above Theorem 
\ref{lem:shape};
\item the diameter of the TD pair $x_{01}, x_{23}$ on $V$, in the
sense of Section 1.
\end{enumerate}
\end{lemma}
\noindent {\it Proof:} 
In the context of either (i), (ii) above
the diameter is one less than the number 
of distinct eigenvalues for $x_{01}$ on $V$.
\hfill $\Box $ \\

\begin{lemma}
\label{lem:double}
Let $V$ denote a finite-dimensional irreducible
$\boxtimes$-module. Then the following coincide:
\begin{enumerate}
\item the shape of the $\boxtimes$-module $V$ from above Theorem 
\ref{lem:shape};
\item the shape of the TD pair $x_{01}, x_{23}$ on $V$, in the
sense of Section 1.
\end{enumerate}
\end{lemma}
\noindent {\it Proof:} 
With reference to 
Lemma \ref{lem:diam2}
let $d$ denote the diameter.
In the context of either (i), (ii) above,
for $0 \leq n \leq d$ the $n$th component of the shape
vector is equal to the dimension of the eigenspace
for $x_{01}$ on $V$ associated with the eigenvalue $2n-d$.
\hfill $\Box $ \\

\begin{theorem}
\label{thm:shape}
Let $\lbrace \rho_i\rbrace_{i=0}^d$ denote the shape
of a tridiagonal pair over $\K$
 that has Krawtchouk type.
 Then  there exists a nonnegative integer $N$ and
  positive integers $d_1, d_2, \ldots, d_N$ such that
\begin{eqnarray*}
  \sum_{i=0}^d \rho_i \lambda^i
  = \prod_{j=1}^N
   (1+\lambda+\lambda^2+\cdots + \lambda^{d_j}).
   \end{eqnarray*}
In particular $\rho_0=1$.
 \end{theorem}
\noindent {\it Proof:} 
Let $A,A^*$ denote a TD pair over $\K$  that has
shape
$\lbrace \rho_i\rbrace_{i=0}^d$ 
and let $V$ denote the underlying vector space.
By Theorem
\ref{prop:ha2}
there exists an irreducible $\boxtimes$-module 
structure on $V$ such that the generators
$x_{01}, x_{23}$ act as
$A, A^*$ respectively.
By Lemma
\ref{lem:double}
the shape of the $\boxtimes$-module $V$
is equal to 
$\lbrace \rho_i\rbrace_{i=0}^d$.
The result follows in view of
Theorem
\ref{lem:shape}.
\hfill $\Box $ \\

\section{TD pairs of Krawtchouk type; the bilinear form}

\noindent In this section we describe a bilinear
form associated with a TD pair of Krawtchouk type.

\begin{theorem}
\label{thm:bilmain}
Let $V$ denote a vector space over $\K$ with finite
positive dimension and
let $A,A^*$ denote a TD pair on $V$ that has Krawtchouk type. Then
there exists a nonzero
bilinear form $\langle\,,\,\rangle$ on $V$ such that both
\begin{eqnarray}
\langle Au,v\rangle= \langle u,Av \rangle,
\qquad \qquad
\langle A^*u,v \rangle= \langle u,A^*v \rangle
\label{eq:dual}
\end{eqnarray}
for $u,v\in V$.
This form is unique up to multiplication by a nonzero scalar in 
$\F$. This form is nondegenerate and symmetric.
\end{theorem}
\noindent {\it Proof:} 
By Theorem \ref{prop:ha2}
there exists a $\boxtimes$-module structure
on $V$ such that the generators $x_{01}, x_{23}$ act as $A, A^*$
respectively.
Consider the element $\sigma = (01)(23)$ in $G$
and let 
$\langle \,,\,\rangle =
\langle \,,\,\rangle_{\sigma}$
denote the corresponding bilinear form
on $V$ from
Theorem \ref{thm:bill}.
We show 
$\langle \,,\,\rangle$ satisfies 
(\ref{eq:dual}).
Observe 
\begin{eqnarray*}
\langle Au,v\rangle &=& \langle x_{01}.u,v\rangle \\
 &=& -\langle u,\sigma(x_{01}).v\rangle \\
 &=& -\langle u,x_{10}.v\rangle \\
 &=& \langle u,x_{01}.v\rangle \\
 &=& \langle u,Av\rangle 
\end{eqnarray*}
and similarly
\begin{eqnarray*}
\langle A^*u,v\rangle &=& \langle x_{23}.u,v\rangle \\
 &=& -\langle u,\sigma(x_{23}).v\rangle \\
 &=& -\langle u,x_{32}.v\rangle \\
 &=& \langle u,x_{23}.v\rangle \\
 &=& \langle u,A^*v\rangle. 
\end{eqnarray*}
Therefore $\langle \,,\,\rangle $
satisfies (\ref{eq:dual}).
By Theorem \ref{thm:bill} the form
$\langle \,,\,\rangle$ is nondegenerate and symmetric.
Concerning the uniqueness of
$\langle \,,\,\rangle $
let $\langle \,,\,\rangle'$
denote any bilinear form
on $V$ that satisfies
(\ref{eq:dual}).
We show that
$\langle \,,\,\rangle'$ is a scalar multiple
of 
$\langle \,,\,\rangle$.
Pick a basis for $V$,
and  let $A_b$ (resp. $A^*_b$) denote the 
matrix that represents $A$ (resp. $A^*$)
with respect to
this basis.
Let $M$ (resp. $N$) denote the matrix that represents 
$\langle \,,\,\rangle$  
(resp. 
$\langle \,,\,\rangle'$)  
with respect to the basis.
Note that $M$ is invertible since 
$\langle \,,\,\rangle$  is nondegenerate.
By 
(\ref{eq:dual})
we have 
$\xi^t_b M = M\xi_b$
and  $\xi^t_b N = N\xi_b$ for $\xi \in 
\lbrace A,A^*\rbrace$.
Combining these equations we find that
$M^{-1}N$ 
commutes with each of
$A_b, A^*_b$.
Let $\gamma $ denote the element of 
${\rm End}(V)$ that is represented by $M^{-1}N$
with respect to the above basis. Then
$\gamma$ commutes with each of $A,A^*$ and
is therefore a scalar multiple of the identity
by Lemma
\ref{lem:id}.
Now $M^{-1}N$ is a scalar multiple of the identity
so $N$ is a scalar multiple of $M$
and therefore 
$\langle \,,\,\rangle'$  is a scalar multiple of
$\langle \,,\,\rangle$.
It follows that
$\langle \,,\,\rangle$ is unique up to multiplication
by a nonzero scalar in $\K$.
\hfill $\Box $ \\

\section{TD pairs of Krawtchouk type; the map $\dagger$}

\noindent In this section we describe an antiautomorphism $\dagger$
associated with a TD pair of Krawtchouk type. We start with a definition.

\begin{definition} 
\label{def:dagger}
Let $V$ denote a vector space over $\K$ with 
finite positive dimension and 
let $A,A^*$ denote a TD pair on $V$ that has
Krawtchouk type. Let $\dagger$ denote the antiautomorphism
of 
${\rm End}(V)$ that corresponds to the bilinear form
in Theorem
\ref{thm:bilmain}. 
We emphasize
\begin{eqnarray}
\langle Xu,v\rangle = 
\langle u,X^\dagger v\rangle 
\qquad \qquad X \in {\rm End}(V), \quad u,v \in V.
\label{eq:dag}
\end{eqnarray}
\end{definition}

\begin{theorem}
\label{thm:dag}
Let $V$ denote a vector space over $\K$ with
finite positive dimension and 
let $A,A^*$ denote a TD pair on $V$ that has Krawtchouk type.
Then $\dagger$ is the unique  antiautomorphism of ${\rm End}(V)$
that fixes each of $A,A^*$. Moreover $X^{\dagger \dagger}=X $
for all $X \in {\rm End}(V)$.
\end{theorem}
\noindent {\it Proof:} 
We first show that $A^\dagger=A$.
For $u,v \in V$ we have
$\langle Au,v\rangle = \langle u,Av\rangle $ by
(\ref{eq:dual}) and
$\langle Au,v\rangle = \langle u,A^\dagger v\rangle $ by
(\ref{eq:dag}) so
 $\langle u,(A^\dagger-A)v\rangle =0$. 
Now  
 $(A^\dagger-A)v =0$ since  
 $\langle \,,\,\rangle $ is nondegenerate
and therefore $A^\dagger=A$.
Similarly
 $A^{*\dagger}=A^*$.
Concerning the uniqueness of $\dagger$,
let $\psi$ denote any antiautomorphism of
${\rm End}(V)$ that fixes each of $A,A^*$.
We show $\psi=\dagger$. The composition
$\psi^{-1} \dagger$ is an automorphism
of ${\rm End}(V)$ that fixes each of $A,A^*$.
This composition is the
identity by Corollary
\ref{cor:id2} so
$\psi=\dagger$. Also since
$\dagger^{-1}$ is an antiautomorphism
of 
${\rm End}(V)$ that fixes each of $A,A^*$
we can take
$\psi=\dagger^{-1}$ in our above argument
to get 
$\dagger^{-1}=\dagger$. This gives
$X^{\dagger \dagger}=X $
for all $X \in {\rm End}(V)$.
\hfill $\Box $ \\ 

\section{TD pairs of Krawtchouk type; the dual space}

\noindent In this section we show that each TD pair of
Krawtchouk type is isomorphic to its dual TD pair in the
sense of 
\cite[Theorem~1.2]{CurtH}.

\begin{theorem} 
\label{lem:chas} Let $V$ denote a vector space over $\K$
with finite positive dimension and let $A,A^*$ denote
a TD pair on $V$ that has Krawtchouk type.
Let $B$ (resp. $B^*$) denote the image
of $A$ (resp. $A^*$) under the canonical 
$\K$-algebra anti-isomorphism
${\rm End}(V)\to 
{\rm End}({\tilde V})$, where
${\tilde V}$ is the vector space dual to $V$.
Then the TD pairs
$A,A^*$ and 
 $B,B^*$ are isomorphic.
\end{theorem}
\noindent {\it Proof:} 
By
 Theorem
\ref{thm:dag} the map
$\dagger$ is an antiautomorphism
of ${\rm End}(V)$ that fixes each of
$A,A^*$. By construction the canonical $\K$-algebra
anti-isomorphism
 ${\rm End}(V)\to {\rm End}({\tilde V})$
sends $A$ to $B$ and $A^*$ to $B^*$.
The
 composition of these two maps
is an $\K$-algebra isomorphism
${\rm End}(V)
\to {\rm End}({\tilde V})$
that sends
$A$ to $B$ and
$A^*$ to  
 $B^*$.
The result follows in view of
Corollary 
\ref{cor:iso2}.
\hfill $\Box $ \\

\section{TD pairs of Krawtchouk type; some isomorphisms}

\noindent In this section we display some isomorphisms
between TD pairs of Krawtchouk type.

\begin{theorem}
\label{thm:4iso}
Let $A,A^*$ denote a TD pair over $\K$ that has Krawtchouk type.
Then the following TD pairs are mutually isomorphic:
\begin{eqnarray*}
A,A^*; \qquad \qquad 
-A,-A^*; \qquad \qquad 
A^*,A; \qquad \qquad 
-A^*,-A.
\end{eqnarray*}
\end{theorem}
\noindent {\it Proof:} 
We first show that the TD pairs $A,A^*$ and $-A,-A^*$ are
isomorphic. Let $V$ denote the vector space underlying
$A,A^*$.
By 
Theorem
\ref{prop:ha2}
there exists a $\boxtimes$-module structure
on $V$ such that the generators $x_{01}, x_{23}$ act as $A, A^*$
respectively.
 Consider the element
$\sigma = (01)(23)$ in $G$.
By 
Theorem \ref{thm:G}
the $\boxtimes$-module $V$
is isomorphic to the $\boxtimes$-module $V$ twisted via $\sigma$.
Let $\gamma:V\to V$ denote
an isomorphism of $\boxtimes$-modules from
$V$ to 
$V$
twisted via $\sigma$.
By the definition of twisting and since $\sigma^2=1$,
\begin{eqnarray*}
\gamma \xi.v = \sigma(\xi).\gamma v \qquad \qquad 
\xi \in \boxtimes, \quad v \in V.
\end{eqnarray*}
We show that $\gamma$ is an isomorphism of TD pairs from
$A,A^*$ to $-A,-A^*$. By construction $\gamma$ is
an isomorphism of vector spaces from $V$ to $V$.
For $v \in V$,
\begin{eqnarray*}
\gamma Av &=& \gamma x_{01}.v \\
           &=& \sigma(x_{01}).\gamma v \\
           &=& x_{10}.\gamma v \\
           &=& -x_{01}.\gamma v \\
           &=& - A\gamma v
\end{eqnarray*}
so $\gamma A=-A \gamma$.
Similarly
\begin{eqnarray*}
\gamma A^*v &=& \gamma x_{23}.v \\
           &=& \sigma(x_{23}).\gamma v \\
           &=& x_{32}.\gamma v \\
           &=& -x_{23}.\gamma v \\
           &=& - A^*\gamma v
\end{eqnarray*}
so $\gamma A^*=-A^* \gamma$.
By the above comments and Definition
\ref{def:iso} the map
$\gamma$ is an isomorphism of
TD pairs from $A,A^*$ to $-A,-A^*$. Therefore the TD pairs
$A,A^*$ and $-A,-A^*$ are isomorphic.
The above argument with $\sigma=(02)(13)$
(resp. $\sigma=(03)(12)$) 
shows that the TD pairs $A,A^*$ and $A^*,A$ (resp.
 $A,A^*$ and $-A^*,-A$)
are isomorphic.
\hfill $\Box $ \\

\section{TD pairs of Krawtchouk type; the Drinfel'd polynomial}

\noindent In this section we introduce a Drinfel'd polynomial
for TD pairs of Krawtchouk type.

\begin{definition}
\label{def:drinv2}
\rm 
Let $A,A^*$ denote a TD pair over $\K$ that has Krawtchouk type.
We define a polynomial $P=P_{A,A^*}$ in $\K\lbrack \lambda \rbrack$
by
\begin{eqnarray*}
P= \sum_{i=0}^d \frac{(-1)^i \zeta_i \lambda^i}{(i!)^2 4^i}
\end{eqnarray*}
where 
$\lbrace \zeta_i \rbrace_{i=0}^d$ is the
split sequence for $A,A^*$ associated with the
standard ordering $\lbrace d-2i \rbrace_{i=0}^d$
(resp. 
 $\lbrace 2i-d\rbrace_{i=0}^d$)
of the eigenvalues for $A$ (resp. $A^*$).
We call $P$ the {\it Drinfel'd polynomial} of $A,A^*$. 
\end{definition}

\begin{lemma}
\label{lem:drin}
Let $V$ denote a finite-dimensional irreducible
$\boxtimes$-module.
Then the following coincide:
\begin{enumerate}
\item the Drinfel'd polynomial $P_V$ from 
line (\ref{eq:drin});
\item the Drinfel'd polynomial for the TD pair $x_{01}, x_{23}$ on $V$,
as in 
Definition
\ref{def:drinv2}.
\end{enumerate}
\end{lemma}
\noindent {\it Proof:}
Let $A:V\to V$ (resp. $A^*:V\to V$) denote the action
of $x_{01}$ (resp. $x_{23}$) on $V$. 
By Theorem
\ref{prop:ha1}
the pair
$A,A^*$ is a TD pair on $V$ that has Krawtchouk type.
We show $P_V=P_{A,A^*}$.
Let $d$ denote the diameter of $A,A^*$
and abbreviate $\theta_i=d-2i$ and
$\theta^*_i=2i-d$ for $0 \leq i \leq d$.
Let $\lbrace U_i\rbrace_{i=0}^d$ denote the split
decomposition  of $V$ for $A,A^*$ and with respect to
$\lbrace \theta_i\rbrace_{i=0}^d$,
$\lbrace \theta^*_i\rbrace_{i=0}^d$.
By \cite[Section~3]{Ev} for $0 \leq i \leq d$
the space $U_i$ is the eigenspace 
for $x_{20}$ on $V$ associated with the eigenvalue
$2i-d$. In particular 
$U_0$ is the eigenspace for $x_{20}$ on $V$ associated with 
the eigenvalue $-d$. Therefore 
$U_0$ coincides with the space $U$
from below
(\ref{eq:epm}).
Using (\ref{eq:epm}) 
we observe that for $0 \leq i \leq d$
the action of $2e^+$ (resp. $2e^-$) on $U_i$ coincides
with the restriction of 
$A-\theta_i I$ (resp. $A^*-\theta^*_i I$) to $U_i$.
By this and 
(\ref{eq:ath}), 
(\ref{eq:asths}) we find
$e^+U_i\subseteq U_{i+1}$ and $e^-U_i\subseteq U_{i-1}$.
By these comments
the action of $4^i(e^-)^i(e^+)^i$ on
$U_0$ coincides with the restriction of
the linear transformation (\ref{eq:aaseig}) to $U_0$. Therefore
$4^i \vartheta_i =
\zeta_i $
where 
$\vartheta_i$ is from
above 
(\ref{eq:drin}) and
$\zeta_i$ is from
below
(\ref{eq:aaseig}).
Comparing
(\ref{eq:drin})
with Definition \ref{def:drinv2} we find
$P_V=P_{A,A^*}$ and the result follows.
\hfill $\Box $ \\

\begin{theorem}
\label{thm:drinf}
The map $A,A^* \mapsto P_{A,A^*}$ induces a bijection between
the following two sets:
\begin{enumerate}
\item the isomorphism classes of TD pairs over $\F$ that have
Krawtchouk type;
\item the polynomials in $\K\lbrack \lambda \rbrack$ that
have constant coefficient 1 and are nonzero at $\lambda=1$.
\end{enumerate}
\end{theorem}
\noindent {\it Proof:} 
The composition of the bijection in
Remark
\ref{prop:ha3} and the bijection in
Theorem
\ref{thm:drintet} is a bijection
from set (i) to set (ii) above.
The map $A,A^* \mapsto P_{A,A^*}$ induces this bijection 
in view of Lemma
\ref{lem:drin}.
\hfill $\Box $ \\

\section{Directions for further research}

\noindent In Sections 8--13 we obtained some results
that apply to TD pairs of Krawtchouk type.
In this section we consider whether these results apply
to more general TD pairs. We also consider
some formulae involving the split sequence.
We content ourselves with
a list of conjectures and problems.

\medskip
\noindent Throughout this section the field $\K$ is arbitrary.

\begin{conjecture}
\label{conj:shape}
Let $\lbrace \rho_i\rbrace_{i=0}^d$ denote the shape
of a tridiagonal pair over $\K$.
 Then  there exists a nonnegative integer $N$ and
  positive integers $d_1, d_2, \ldots, d_N$ such that
\begin{eqnarray*}
  \sum_{i=0}^d \rho_i \lambda^i
  = \rho_0 \prod_{j=1}^N
   (1+\lambda+\lambda^2+\cdots + \lambda^{d_j}).
   \end{eqnarray*}
Moreover if $\K$ is algebraically closed then
$\rho_0=1$.
\end{conjecture}

\begin{conjecture}
\label{conj:bilmain}
Assume $\K$ is algebraically closed.
Let $V$ denote a vector space over $\K$ with finite
positive dimension and
let $A,A^*$ denote a TD pair on $V$.
 Then
there exists a nonzero
bilinear form $\langle\,,\,\rangle$ on $V$ such that both
\begin{eqnarray*}
\langle Au,v\rangle= \langle u,Av \rangle,
\qquad \qquad
\langle A^*u,v \rangle= \langle u,A^*v \rangle
\end{eqnarray*}
for $u,v\in V$.
This form is unique up to multiplication by a nonzero scalar in 
$\F$. This form is nondegenerate and symmetric.
\end{conjecture}

\begin{conjecture}
\label{conj:dag}
Assume $\K$ is algebraically closed.
Let $V$ denote a vector space over $\K$ with
finite positive dimension and 
let $A,A^*$ denote a TD pair on $V$.
Then there exists a  unique  antiautomorphism $\dagger$ of ${\rm End}(V)$
that fixes each of $A,A^*$. Moreover $X^{\dagger \dagger}=X $
for all $X \in {\rm End}(V)$.
\end{conjecture}

\begin{conjecture} 
\label{conj:chas}
Assume $\K$ is algebraically closed.
Let $V$ denote a vector space over $\K$
with finite positive dimension and let $A,A^*$ denote
a TD pair on $V$.
Let $B$ (resp. $B^*$) denote the image
of $A$ (resp. $A^*$) under the canonical 
$\K$-algebra anti-isomorphism
${\rm End}(V)\to 
{\rm End}({\tilde V})$,
where ${\tilde V}$ is the vector space dual to $V$.
Then the TD pairs
$A,A^*$ and 
 $B,B^*$ are isomorphic.
\end{conjecture}

\noindent 
In order to state the next conjecture
we make a definition.
\begin{definition}
\label{def:pararray}
\rm
Let $A,A^*$ denote a TD pair and
assume the shape vector satisfies
$\rho_0=1$.
By a {\it parameter array} of $A,A^*$ we mean
a sequence 
$(\lbrace \theta_i\rbrace_{i=0}^d;
\lbrace \theta^*_i\rbrace_{i=0}^d;
\lbrace \zeta_i\rbrace_{i=0}^d)$
where 
$\lbrace \theta_i\rbrace_{i=0}^d$ (resp.
$\lbrace \theta^*_i\rbrace_{i=0}^d$)
is a standard ordering of the eigenvalues of
$A$ (resp. $A^*$) and
$\lbrace \zeta_i\rbrace_{i=0}^d$ is the
split sequence of $A,A^*$
with respect  to
$\lbrace \theta_i\rbrace_{i=0}^d$
and $\lbrace \theta^*_i\rbrace_{i=0}^d$.
\end{definition}

\begin{conjecture}
\label{conj:main}
Assume $\K$ is algebraically closed.
Let $d$ denote a nonnegative integer and let
\begin{eqnarray}
\label{eq:pa}
(\lbrace \theta_i\rbrace_{i=0}^d;
\lbrace \theta^*_i\rbrace_{i=0}^d;
\lbrace \zeta_i\rbrace_{i=0}^d)
\end{eqnarray}
denote a sequence of scalars taken from $\K$.
Then there exists a TD pair $A,A^*$ over $\K$ with
parameter array 
(\ref{eq:pa}) if and only if (i)--(iii) hold below.
\begin{enumerate}
\item $\zeta_0=1$, $\zeta_d\not=0$, and
\begin{eqnarray*}
0 \not=\sum_{i=0}^d
\zeta_i(\theta_0-\theta_{i+1})\cdots 
(\theta_0-\theta_d)
(\theta^*_0-\theta^*_{i+1})\cdots 
(\theta^*_0-\theta^*_d);
\end{eqnarray*}
\item $\theta_i\not=\theta_j$, 
          $\theta^*_i\not=\theta^*_j$ if $i\not=j$ $(0 \leq i,j\leq d)$;
\item the expressions
\begin{eqnarray*}
\frac{\theta_{i-2}-\theta_{i+1}}{\theta_{i-1}-\theta_i},
\qquad \qquad 
\frac{\theta^*_{i-2}-\theta^*_{i+1}}{\theta^*_{i-1}-\theta^*_i}
\end{eqnarray*}
are equal and independent of $i$ for $2 \leq i \leq d-1$.
\end{enumerate}
Suppose (i)--(iii) hold. Then $A,A^*$ is unique up to isomorphism
of TD pairs.
\end{conjecture}

\noindent As a first step in the proof of
Conjecture
\ref{conj:main}, try to prove the following conjecture.

\begin{conjecture}
Assume $\K$ is algebraically closed. Then two
TD pairs over $\K$ are isomorphic if and only if 
they have a parameter array in common.
\end{conjecture}

\begin{conjecture}
\label{conjz}
Let $A,A^*$ denote a TD pair.
Assume the shape vector satisfies
$\rho_0=1$
and let 
$(\lbrace \theta_i\rbrace_{i=0}^d;
\lbrace \theta^*_i\rbrace_{i=0}^d;
\lbrace \zeta_i\rbrace_{i=0}^d)$ denote 
a parameter array of $A,A^*$.
Then for $0 \leq i \leq d$ both 
\begin{eqnarray*}
\zeta_i &=& (\theta^*_0-\theta^*_1)
(\theta^*_0-\theta^*_2)\cdots 
(\theta^*_0-\theta^*_i)
{\rm tr}(\tau_i(A)E^*_0),
\\
\zeta_i &=& (\theta_0-\theta_1)
(\theta_0-\theta_2)\cdots 
(\theta_0-\theta_i)
{\rm tr}(\tau^*_i(A^*)E_0)
\end{eqnarray*}
where $E_0$ (resp. $E^*_0$) is the primitive idempotent
of $A$ (resp. $A^*$) for $\theta_0$ (resp. $\theta^*_0$)
\cite[Section~2]{TD00}
and
\begin{eqnarray*}
\tau_i &=& (\lambda-\theta_0)
(\lambda-\theta_1) \cdots
(\lambda-\theta_{i-1}),
\\
\tau^*_i &=& (\lambda-\theta^*_0)
(\lambda-\theta^*_1) \cdots
(\lambda-\theta^*_{i-1}).
\end{eqnarray*}
\end{conjecture}

\begin{conjecture}
With the notation and assumptions of
Conjecture
\ref{conjz} we have 
 ${\rm {tr}}(E_0E^*_0)\not=0$. Moreover $\zeta_i$ is equal to
 each of
\begin{eqnarray*}
 \frac{{\rm {tr}}(E_0\tau^*_i(A^*)\tau_i(A)E^*_0)}{{\rm {tr}}(E_0E^*_0)},
\qquad \qquad 
 \frac{{\rm {tr}}(E^*_0\tau_i(A)\tau^*_i(A^*)E_0)}{{\rm {tr}}(E^*_0E_0)}
\end{eqnarray*}
for $0 \leq i \leq d$.
\end{conjecture}

\begin{problem}
\rm
\label{probzz}
Let $A,A^*$ denote a TD pair.
Let $\lbrace \theta_i\rbrace_{i=0}^d$
(resp. 
$\lbrace \theta^*_i\rbrace_{i=0}^d$)
denote a standard ordering of the eigenvalues
of $A$ (resp. $A^*$). 
Assume the shape vector satisfies 
$\rho_0=1$.
Find the algebraic relationships between the following
eight sequences.
\begin{enumerate}
\item
The split sequences of $A,A^*$ with respect to
\begin{eqnarray*}
&&\lbrace \theta_i\rbrace_{i=0}^d
\; {\mbox{\rm and}} \; \lbrace \theta^*_i\rbrace_{i=0}^d;
\qquad \qquad \;\;
\lbrace \theta_{d-i}\rbrace_{i=0}^d
\; {\mbox{\rm and}} \; \lbrace \theta^*_i\rbrace_{i=0}^d;
\\
&&\lbrace \theta_i\rbrace_{i=0}^d
\; {\mbox{\rm and}} \; \lbrace \theta^*_{d-i}\rbrace_{i=0}^d;
\qquad \qquad 
\lbrace \theta_{d-i}\rbrace_{i=0}^d
\; {\mbox{\rm and}} \; \lbrace \theta^*_{d-i}\rbrace_{i=0}^d;
\end{eqnarray*}
\item
The split sequences of $A^*,A$ with respect to
\begin{eqnarray*}
&&\lbrace \theta^*_i\rbrace_{i=0}^d
\; {\mbox{\rm and}} \; \lbrace \theta_i\rbrace_{i=0}^d;
\qquad \qquad \;\; \;
\lbrace \theta^*_{d-i}\rbrace_{i=0}^d
\; {\mbox{\rm and}} \; \lbrace \theta_i\rbrace_{i=0}^d;
\\
&&\lbrace \theta^*_i\rbrace_{i=0}^d
\; {\mbox{\rm and}}\; \lbrace \theta_{d-i}\rbrace_{i=0}^d;
\qquad \qquad 
\lbrace \theta^*_{d-i}\rbrace_{i=0}^d
\; {\mbox{\rm and}} \; \lbrace \theta_{d-i}\rbrace_{i=0}^d.
\end{eqnarray*}
\end{enumerate}
\end{problem}

\begin{conjecture}
\label{conjzz}
Let $A,A^*$ denote a TD pair.
Let 
$\lbrace \theta_i\rbrace_{i=0}^d$
(resp. 
$\lbrace \theta^*_i\rbrace_{i=0}^d$)
denote a standard ordering of the eigenvalues
of $A$ (resp. $A^*$). 
Assume the shape vector satisfies 
$\rho_0=1$.
Then the following coincide:
\begin{enumerate}
\item
The split sequence of $A,A^*$ with respect to
$\lbrace \theta_i\rbrace_{i=0}^d$
and $\lbrace \theta^*_i\rbrace_{i=0}^d$;
\item
The split sequence of $A^*,A$ with respect to
$\lbrace \theta^*_i\rbrace_{i=0}^d$
and $\lbrace \theta_i\rbrace_{i=0}^d$.
\end{enumerate}
\end{conjecture}

\begin{problem}
\rm
Referring to Proposition 
\ref{prop:id} and Corollary
\ref{cor:id2}, what conclusions can we
obtain if we drop the assumption that
$\K$ is algebraically closed?
\end{problem}

\begin{problem} \rm
Referring to Theorem
\ref{thm:4iso}, 
let $B,B^*$ denote one of the TD pairs listed
in that theorem, other than $A,A^*$. 
Let $\gamma $ denote an isomorphism of  
TD pairs from $A,A^*$ to $B,B^*$.
Express $\gamma$ explicitly as a polynomial in $A,A^*$. 
\end{problem}

\noindent Tatsuro Ito \hfil\break
\noindent Department of Computational Science \hfil\break
\noindent Faculty of Science \hfil\break
\noindent Kanazawa University \hfil\break
\noindent Kakuma-machi \hfil\break
\noindent Kanazawa 920-1192, Japan \hfil\break
\noindent email:  {\tt ito@kappa.s.kanazawa-u.ac.jp}

\bigskip

\noindent Paul Terwilliger \hfil\break
\noindent Department of Mathematics \hfil\break
\noindent University of Wisconsin \hfil\break
\noindent 480 Lincoln Drive \hfil\break
\noindent Madison, WI 53706-1388 USA \hfil\break
\noindent email: {\tt terwilli@math.wisc.edu }\hfil\break


\begin{thebibliography}{10}


%

\bibitem{hasan}
H.~Alnajjar  and B.~Curtin.
\newblock
A family of tridiagonal pairs.
\newblock {\em
Linear Algebra Appl.}
{\bf 390}
(2004)
369--384.

\bibitem{hasan2}
H.~Alnajjar  and B.~Curtin.
\newblock
A family of tridiagonal pairs related to
the quantum affine algebra 
$U\sb q(\widehat{\mathfrak{sl}}\sb 2)$.
\newblock {\em
Electron. J. Linear Algebra}
{\bf 13} (2005) 1--9. 

\bibitem{CurtH}
H.~Alnajjar  and B.~Curtin.
\newblock
A bilinear form for tridiagonal pairs of $q$-Serre type.
\newblock {\em
Linear Algebra Appl.}, submitted.




\bibitem{AWil}
R.~Askey and J.A.~Wilson.
\newblock A set of orthogonal polynomials that 
generalize the {R}acah coefficients or $6-j$ symbols. 
\newblock {\em SIAM J. Math. Anal.}, 10:1008--1016, 1979. 



							       
							    





							       





%




\bibitem{BT}
G.~Benkart and P.~Terwilliger.
\newblock The universal central extension of the three-point
$\mathfrak{sl}_2$ loop algebra.
\newblock{\em Proc. Amer. Math. Soc.} {\bf 135} (2007) 1651--1657. 
{\tt arXiv:math.RA/0512422}.

%





























\bibitem{E}
A.~Elduque.
\newblock The $S_4$-action on the tetrahedron algebra.
\newblock Preprint. \hfil\break
{\tt arXiv:math.RA/0604218}.





\bibitem{GR}
G.~Gasper and M.~Rahman.
\newblock {\em Basic hypergeometric series}.
\newblock Encyclopedia of Mathematics and its Applications, 35,
\newblock Cambridge University Press, Cambridge, 1990.



















\bibitem{Ha}
B.~Hartwig.
\newblock The tetrahedron algebra and its finite-dimensional
irreducible modules. 
{\em Linear Algebra Appl.} {\bf 422} (2007) 219--235;
{\tt arXiv:math.RT/0606197}.

\bibitem{HT}
B.~Hartwig and P.~Terwilliger.
\newblock The Tetrahedron algebra, the Onsager algebra,
and the $\mathfrak{sl}_2$ loop algebra.
\newblock{\em J. Algebra} {\bf 308} (2007) 840--863;
{\tt arXiv:math.ph/0511004}.






\bibitem{TD00}
T.~Ito, K.~Tanabe, and P.~Terwilliger.
\newblock Some algebra related to ${P}$- and ${Q}$-polynomial association
  schemes,  in:
\newblock {\em Codes and Association Schemes (Piscataway NJ, 1999)}, Amer.
Math. Soc., Providence RI, 2001, pp.
     167--192; 
{\tt arXiv:math.CO/0406556}.

\bibitem{shape}
T.~Ito and P.~Terwilliger.
\newblock The shape of a tridiagonal pair.
\newblock {\em J. Pure Appl. Algebra}
	      {\bf 188}
	    (2004)
		     145--160;
{\tt arXiv:math.QA/0304244
}.

\bibitem{tdanduq}
T.~Ito and P.~Terwilliger.
\newblock {Tridiagonal pairs and the quantum affine 
algebra
$U_q({\widehat{sl}}_2)$.}
\newblock {\em Ramanujan J.} 
{\bf 13} (2007) 39--62;
{\tt arXiv:math.QA/0310042}.

\bibitem{NN}
T.~Ito and P.~Terwilliger.
\newblock Two non-nilpotent linear transformations that
satisfy the cubic $q$-Serre relations.
\newblock{\em J. Algebra Appl.}, submitted;
{\tt arXiv:math.QA/050839}.

\bibitem{qtet}
T.~Ito and P.~Terwilliger.
\newblock The $q$-tetrahedron algebra and its
finite-dimensional irreducible modules.
\newblock{\em Comm. Algebra}, to appear; 
{\tt arXiv:math.QA/0602199}.

\bibitem{ITdrg}
T.~Ito and P.~Terwilliger.
\newblock Distance-regular graphs and the 
$q$-tetrahedron algebra.
\newblock{\em European J. Combin.}, submitted;
{\tt arXiv:math.CO/0608694}.


\bibitem{Ev}
T.~Ito and P.~Terwilliger.
\newblock Finite-dimensional irreducible modules for
the three-point $\mathfrak{sl}_2$ loop algebra.
Preprint.







						  
\bibitem{KoeSwa}
R.~Koekoek and R.~F.~Swarttouw.
\newblock {\em The Askey scheme of hypergeometric orthogonal
polyomials and its
  $q$-analog}, report 98-17, Delft University of Technology, The
  Netherlands, 1998.
  Available at
 \newblock{ \tt http://aw.twi.tudelft.nl/{\~{}}koekoek/research.html} 




















\bibitem{N:aw}
K.~Nomura.
\newblock
Tridiagonal pairs and the {A}skey-{W}ilson relations.
\newblock{\em
Linear Algebra Appl.} {\bf 397} (2005) 99--106.

\bibitem{N:refine}
K.~Nomura.
\newblock
A refinement of the split decomposition of a tridiagonal pair.
\newblock{\em 
Linear Algebra Appl.} {\bf 403} (2005) 1--23.

\bibitem{NT:balanced}
K.~Nomura, P.~Terwilliger.
\newblock
Balanced Leonard pairs.
\newblock{\em 
Linear Algebra Appl.} {\bf 420} (2007) 51--69;
{\tt arXiv:math.RA/0506219}.

\bibitem{NT:formula}
K.~Nomura, P.~Terwilliger.
\newblock
Some trace formulae involving the split sequences of a Leonard pair.
\newblock{\em 
Linear Algebra Appl.} {\bf 413} (2006) 189--201;
{\tt arXiv:math.RA/0508407}.

\bibitem{NT:det}
K.~Nomura, P.~Terwilliger.
\newblock
The determinant of $AA^*-A^*A$ for a Leonard pair $A,A^*$.
\newblock{\em 
Linear Algebra Appl.} {\bf 416} (2006) 880--889;
{\tt arXiv:math.RA/0511641}.

\bibitem{NT:mu}
K.~Nomura, P.~Terwilliger.
\newblock
Matrix units associated with the split basis of a Leonard pair.
\newblock{\em 
Linear Algebra Appl.} {\bf 418} (2006) 775--787;
{\tt arXiv:math.RA/0602416}.

\bibitem{NT:span}
K.~Nomura, P.~Terwilliger.
\newblock
Linear transformations that are tridiagonal with respect to
both eigenbases of a Leonard pair.
\newblock{\em
Linear Algebra Appl.} {\bf 420} (2007) 198--207.
{\tt arXiv:math.RA/0605316}.

\bibitem{NT:switch}
K.~Nomura, P.~Terwilliger.
\newblock
The switching element for a Leonard pair.
\newblock{\em
Linear Algebra Appl.}, submitted for publication;
{\tt arXiv:math.RA/0608623}.

\bibitem{nomsplit}
K.~Nomura and P.~Terwilliger.
\newblock
The split decomposition of a tridiagonal pair.
\newblock {\em Linear Algebra Appl.},
to appear;
{\tt arXiv:math.RA/0612460}. 









\bibitem{PT}
A.~A.~Pascasio and P.~Terwilliger.
\newblock The tetrahedron algebra and the Hamming graphs.
\newblock In preparation.







\bibitem{rotman}
J.~J.~Rotman.
\newblock{\em Advanced modern algebra.}
\newblock{Prentice Hall, Saddle River NJ 2002}.










   


\bibitem{LS99}
P.~Terwilliger.
\newblock Two linear transformations each tridiagonal with respect to an
  eigenbasis of the other.
  \newblock {\em Linear Algebra Appl.}  {\bf 330} (2001) 149--203;
{\tt arXiv:math.RA/0406555}.

  \bibitem{qSerre}
  P.~Terwilliger.
  \newblock Two relations that generalize the $q$-Serre relations and the
  Dolan-Grady relations. In
  \newblock {\em  Physics and
  Combinatorics 1999 (Nagoya)}, 377--398, World Scientific Publishing,
   River Edge, NJ, 2001; 
{\tt arXiv:math.QA/0307016}.

   \bibitem{LS24}
   P.~Terwilliger.
   \newblock  Leonard pairs from 24 points of view.
   \newblock {\em Rocky Mountain J. Math.} {\bf 32}(2) (2002) 827--888;
{\tt arXiv:math.RA/0406577}.

   \bibitem{conform}
   P.~Terwilliger.
   \newblock Two linear transformations each tridiagonal with respect to an
     eigenbasis of the other; the $TD$-$D$ and the $LB$-$UB$ canonical form.
\newblock {\em J. Algebra} {\bf 298} (2006) 302--319.
{\tt arXiv:math.RA/0304077}.


    \bibitem{lsint}
    P.~Terwilliger.
    \newblock Introduction to Leonard pairs.
    \newblock {OPSFA Rome 2001}.
    \newblock{\em J. Comput. Appl. Math.} {\bf 153}(2) (2003)
    463--475.


\bibitem{TLT:split}
P.~Terwilliger.
\newblock Two linear transformations each tridiagonal with respect to an
  eigenbasis of the other; comments on the split decomposition.
\newblock {\em 
 J. Comput. Appl. Math.} {\bf 178} (2005) 437--452;
{\tt arXiv:math.RA/0306290}.

 \bibitem{TLT:array}
 P.~Terwilliger.
 \newblock Two linear transformations each tridiagonal with respect to an
   eigenbasis of the other; comments on the parameter array.
\newblock {\em
Des. Codes Cryptogr.}  {\bf 34}  (2005) 307--332;
{\tt arXiv:math.RA/0306291}.

\bibitem{qrac}
P.~Terwilliger.
\newblock Leonard pairs and the $q$-Racah polynomials.
\newblock {\em Linear Algebra Appl.} {\bf 387} (2004) 235--276;
{\tt arXiv:math.QA/0306301}.




\bibitem{madrid}
P.~Terwilliger.
\newblock
An algebraic approach to the Askey scheme of orthogonal polynomials. 
Orthogonal polynomials and special functions, 
255--330, Lecture Notes in Math., 1883, 
Springer, Berlin, 2006; 
{\tt arXiv:math.QA/0408390}. 


\bibitem{aw}
P.~Terwilliger and R.~Vidunas.
\newblock Leonard pairs and the Askey-Wilson relations.
\newblock {\em J. Algebra Appl.} {\bf 3} (2004) 411--426;
{\tt arXiv:math.QA/0305356}.








%

 \end{thebibliography}
\end{document}